\documentclass[a4paper,fleqn]{cas-dc}

\usepackage[numbers]{natbib}

\usepackage{todonotes}
\usepackage{mathtools}
\usepackage{tikz}
\usetikzlibrary{arrows.meta,positioning}
\usepackage{pgfplots}
\pgfplotsset{compat=1.18}
\usepackage{cleveref}
\crefname{equation}{}{}

\usepackage{xcolor}
\definecolor{color1}{RGB}{0,128,128}
\definecolor{color2}{RGB}{244,224,77}
\definecolor{color3}{RGB}{251,86,7}

\usepackage{contour}
\usepackage{framed} 
\usepackage{multicol} 

\usepackage{nomencl} 
\makenomenclature
\renewcommand*\nompreamble{\begin{multicols}{2}}
\renewcommand*\nompostamble{\end{multicols}}
\renewcommand\nomgroup[1]{%
  \item[\bfseries
  \ifstrequal{#1}{A}{Sets and Indices}{%
  \ifstrequal{#1}{B}{Parameters}{%
  \ifstrequal{#1}{C}{Variables}{}}}%
]}

\usepackage{lipsum}

\usepackage{amsmath}

\newcommand{\FP}{\mathit{F\!P}}
\newcommand{\MU}{\mathit{M\!U}}
\newcommand{\PC}{\mathit{PC}}
\newcommand{\rp}{\mathit{rp}}
\newcommand{\RP}{\mathit{RP}}
\newcommand{\RU}{\mathit{RU}}
\newcommand{\sd}{\mathit{sd}}
\newcommand{\su}{\mathit{su}}

\newcommand{\Markov}{\textsc{Markov Transition}\xspace}
\newcommand{\Cyclic}{\textsc{Cyclic Connection}\xspace}
\newcommand{\NoEnf}{\textsc{No Enforcement}\xspace}
\newcommand{\Reference}{\textsc{Reference Case}\xspace}

\def\tsc#1{\csdef{#1}{\textsc{\lowercase{#1}}\xspace}}
\tsc{WGM}
\tsc{QE}

\begin{document}
\let\WriteBookmarks\relax
\def\floatpagepagefraction{1}
\def\textpagefraction{.001}

\shorttitle{Connecting RPs using Transition Matrices}

\shortauthors{F. C. A. Auer et al.}

\title [mode = title]{Connecting Representative Periods in Energy System Optimization Models using Markov Transition Matrices}  

\tnotemark[1] 

\tnotetext[1]{This work is part of the project iKlimEt (FO999910627), which has received funding in the framework ``Energieforschung'', a research and technology program of the ``Klima- und Energiefonds'' of Austria.}

%
\author[1,2]{Felix C. A. Auer}[orcid=0009-0009-3585-9628]

\cormark[1]


\ead{felix.auer@tugraz.at}


\credit{Conceptualization, Methodology, Software, Writing - Original Draft, Visualization}

\affiliation[1]{organization={Institute of Electricity Economics and Energy Innovation, Graz University of Technology},
            addressline={Inffeldgasse 18}, 
            city={Graz},
            citysep={}, 
            postcode={8010}, 
            country={Austria}}
\affiliation[2]{organization={Research Center ENERGETIC, Graz University of Technology},
            addressline={Rechbauerstraße 12}, 
            city={Graz},
            citysep={}, 
            postcode={8010}, 
            country={Austria}}

\author[1,2]{Robert Gaugl}[orcid=0000-0003-4112-4483]
\ead{robert.gaugl@tugraz.at}
\credit{Conceptualization, Writing - Review \& Editing, Supervision}

\author[1,2]{Thomas Klatzer}[orcid=0000-0001-8041-0496]
\ead{thomas.klatzer@tugraz.at}
\credit{Writing - Review \& Editing, Project administration, Funding acquisition}

\author[3,4]{Diego A. Tejada-Arango}[orcid=0000-0002-3278-9283]

\ead{diego.tejadaarango@tno.nl}

\credit{Conceptualization, Writing - Review \& Editing, Supervision}

\affiliation[3]{organization={Energy \& Materials Transition Unit, Netherlands Organisation for Applied Scientific Research (TNO)},
            addressline={Westerduinweg 3},
            city={Petten},
            citysep={}, 
            postcode={1755 LE},
            country={The Netherlands}}
\affiliation[4]{organization={Instituto de Investigación Tecnológica, Escuela Técnica Superior de Ingeniería, Universidad Pontificia Comillas},
            addressline={Calle del Rey Francisco 4}, 
            city={Madrid},
            citysep={}, 
            postcode={28008}, 
            country={Spain}}

\author[1,2]{Sonja Wogrin}[orcid=0000-0002-3889-7197]
\ead{wogrin@tugraz.at}
\credit{Conceptualization, Writing - Review \& Editing, Supervision, Funding acquisition}

\cortext[1]{Corresponding author}



\begin{abstract}Time series aggregation reduces the computational complexity of large-scale energy system optimization models, but maintaining chronological continuity between the resulting representative periods (RPs) remains a key challenge, as transitions between RPs are typically lost. This causes inaccuracies in storage behavior, unit commitment, and other time-linked aspects of the model. We propose a novel method that uses the Transition Matrix between RPs to link them via probabilistic transitions and expected values. In contrast to existing Transition Matrix approaches that add variables and constraints to reconstruct inter-period chronology (e.g., for seasonal storage), our method reformulates the existing intra-RP constraints at the period boundaries without introducing any additional variables or constraints. It also handles constraints that connect multiple time steps and can be adapted to binary variables. We demonstrate the benefits on an illustrative case study and validate them on the updated IEEE Reliability Test System (RTS-GMLC). The improvement over the state of the art depends on the structure of the Transition Matrix, which can be inspected a priori at no additional data cost. When it is near-diagonal, the established cyclic connection already performs well, whereas for less diagonal matrices the Markov Transition reduces the median operational deviation by up to 80\% (from about 32\% to 6\%). These gains come at practically no extra cost, as the mean computational effort stays below 2\% of the full-model runtime, at most 0.9 percentage points more than the cyclic connection.
\end{abstract}


\begin{keywords}
energy systems\sep transition matrices\sep Markov\sep mixed-integer linear programming\sep optimization\sep time series aggregation
\end{keywords}

\ExplSyntaxOn
\bool_gset_true:N \g_stm_nologo_bool
\ExplSyntaxOff
\maketitle

\settowidth{\nomlabelwidth}{$\rp\in\mathcal{RP}$\hspace{-0.5em}}

\nomenclature[Ag]{$g\in \mathcal{G}$}{Generating units, from 1 to $G$.}
\nomenclature[Ak]{$k\in \mathcal{K}$}{Time steps within each $\rp$, from 1 to $K$.}
\nomenclature[Ak0]{$k_0$}{Predecessor of the first time step $k_1$. Equals $k_K$ (cyclic) or an expected value (Markov).}
\nomenclature[Arp]{$\rp\in \mathcal{RP}$}{Representative periods, from 1 to $\RP$. Typically $\RP\ll T/K$.}
\nomenclature[As]{$s\in \mathcal{S}$}{Storage units, from 1 to $S$.}
\nomenclature[At]{$t\in \mathcal{T}$}{Time steps (in our case in hours), from 1 to $T$.}

\nomenclature[BD]{$\Delta_k$}{Duration of one time step, used for conversion between power and energy (in our case 1 hour).}
\nomenclature[BFP]{$\FP_{\rp,g}$}{Fixed power output of $g$ at the first time step of $\rp$ [MW].}
\nomenclature[BMU]{$\MU_g$}{Minimum Up-Time of $g$ [h].}
\nomenclature[BN]{$N_{\rp',\rp}$}{Number of transitions from $\rp'$ to $\rp$.}
\nomenclature[BP]{$P_{\rp',\rp}$}{Probability that $\rp$ was preceded by $\rp'$, in $[0,1]$.}
\nomenclature[BPC]{$\PC_g$}{Max. production capability of $g$ when committed [MW].}
\nomenclature[BRU]{$\RU_g$}{Ramp-up rate of unit $g$ per time step [MW/h].}

\nomenclature[Cc]{$c_{\rp,k,s}$}{Charging power at $\rp$, $k$ of storage unit~$s$ [MW].}
\nomenclature[Cl]{$l_{\rp,k,s}$}{Storage level at $\rp$, $k$ of storage unit~$s$ [MWh].}
\nomenclature[Cp]{$p_{\rp,k,*}$}{Power output at $\rp$, $k$. Either from a thermal generator~$g$ or from a storage unit~$s$ (i.e., when discharging) [MW].}
\nomenclature[Csd]{$\sd_{\rp,k,g}$}{Shutdown flag at $\rp$, $k$ of unit~$g$ [binary].}
\nomenclature[Csu]{$\su_{\rp,k,g}$}{Start-up flag at $\rp$, $k$ of unit~$g$ [binary].}
\nomenclature[Cu]{$u_{\rp,k,g}$}{Unit-Commitment flag at $\rp$, $k$ of unit~$g$ [binary].}

\begin{table*}[!t]
  \begin{framed}
    \makeatletter\nom@tempdim\nomlabelwidth\relax\makeatother
    \begin{thenomenclature}
    \nomgroup{A}
      \item [{$g\in \mathcal{G}$}]\begingroup Generating units, from 1 to $G$.\nomeqref {0}\nompageref{1}
      \item [{$k\in \mathcal{K}$}]\begingroup Time steps within each $\rp$, from 1 to $K$.\nomeqref {0}\nompageref{1}
      \item [{$k_0$}]\begingroup Predecessor of the first time step $k_1$. Equals $k_K$ (cyclic) or an expected value (Markov).\nomeqref {0}\nompageref{1}
      \item [{$\rp\in \mathcal{RP}$}]\begingroup Representative periods, from 1 to $\RP$. Typically $\RP\ll T/K$.\nomeqref {0}\nompageref{1}
      \item [{$s\in \mathcal{S}$}]\begingroup Storage units, from 1 to $S$.\nomeqref {0}\nompageref{1}
      \item [{$t\in \mathcal{T}$}]\begingroup Time steps (in our case in hours), from 1 to $T$.\nomeqref {0}\nompageref{1}
    \nomgroup{B}
      \item [{$\Delta_k$}]\begingroup Duration of one time step, used for conversion between power and energy (in our case 1 hour).\nomeqref {0}\nompageref{1}
      \item [{$\FP_{\rp,g}$}]\begingroup Fixed power output of $g$ at the first time step of $\rp$ [MW].\nomeqref {0}\nompageref{1}
      \item [{$\MU_g$}]\begingroup Minimum Up-Time of $g$ [h].\nomeqref {0}\nompageref{1}
      \item [{$N_{\rp',\rp}$}]\begingroup Number of transitions from $\rp'$ to $\rp$.\nomeqref {0}\nompageref{1}
      \item [{$P_{\rp',\rp}$}]\begingroup Probability that $\rp$ was preceded by $\rp'$, in $[0,1]$.\nomeqref {0}\nompageref{1}
      \item [{$\PC_g$}]\begingroup Max. production capability of $g$ when committed [MW].\nomeqref {0}\nompageref{1}
      \item [{$\RU_g$}]\begingroup Ramp-up rate of unit $g$ per time step [MW/h].\nomeqref {0}\nompageref{1}
    \nomgroup{C}
      \item [{$c_{\rp,k,s}$}]\begingroup Charging power at $\rp$, $k$ of storage unit~$s$ [MW].\nomeqref {0}\nompageref{1}
      \item [{$l_{\rp,k,s}$}]\begingroup Storage level at $\rp$, $k$ of storage unit~$s$ [MWh].\nomeqref {0}\nompageref{1}
      \item [{$p_{\rp,k,*}$}]\begingroup Power output at $\rp$, $k$. Either from a thermal generator~$g$ or from a storage unit~$s$ (i.e., when discharging) [MW].\nomeqref {0}\nompageref{1}
      \item [{$\sd_{\rp,k,g}$}]\begingroup Shutdown flag at $\rp$, $k$ of unit~$g$ [binary].\nomeqref {0}\nompageref{1}
      \item [{$\su_{\rp,k,g}$}]\begingroup Start-up flag at $\rp$, $k$ of unit~$g$ [binary].\nomeqref {0}\nompageref{1}
      \item [{$u_{\rp,k,g}$}]\begingroup Unit-Commitment flag at $\rp$, $k$ of unit~$g$ [binary].\nomeqref {0}\nompageref{1}
    \end{thenomenclature}
  \end{framed}
\end{table*}

\section{Introduction}
As power systems shift toward high renewable shares, the optimization models needed to plan them are growing in complexity, which strains what is computationally tractable. The European power system exemplifies this, as it is highly interconnected and is subject to climate policies that actively drive the adoption of renewable technologies. The European Green Deal targets a 55\% emission cut by 2030 (compared to 1990) and climate neutrality by 2050~\cite{noauthor_european_nodate}, with energy being responsible for more than 75\% of the EU's greenhouse gas emissions~\cite{noauthor_energy_2022}. Because renewable sources are weather-dependent, their output varies strongly across hours, days, and seasons. Balancing this variability requires more connections between national energy systems and greater reliance on international electricity trading, which pushes models toward a continental scope~\cite{horschPyPSAEurOpenOptimisation2018}. Robust investment decisions further require accounting for long planning horizons and multiple weather years, as results based on only a few years of data are unreliable at high renewable shares~\cite{pfenningerDealingMultipleDecades2017}. Optimizing investment and operation over multiple years at hourly resolution couples thousands of consecutive time steps, and the problem grows further with finer resolution and wider spatial scope, so that models routinely sacrifice temporal or spatial detail to remain tractable~\cite{craigOvercomingDisconnectEnergy2022}. Reducing the temporal dimension while preserving solution accuracy is therefore one central challenge to planning future, highly renewable energy systems.

One possible approach is performing time series aggregation (TSA), which aims to construct aggregated energy system optimization models (ESOMs) with reduced temporal dimensionality that accurately approximate the optimization outcomes of their full-scale counterparts at reduced computational cost~\cite{teichgraeber_time-series_2022}. TSA methods can be split into two categories~\cite{wogrin_time_2023}: \emph{A priori} methods aggregate using only input data (such as demand time series), typically by employing clustering techniques such as k-means, k-medoids, or hierarchical clustering. In contrast, \emph{a posteriori} methods integrate information from the optimization model they aim to approximate into the TSA procedure (e.g., through preliminary optimizations). Despite their demonstrated potential (e.g., in~\cite{wogrin_time_2023}), a posteriori methods have not seen wide adoption in large-scale energy system models, so a priori methods remain the current state of the art.

A priori methods comprise, among others, downsampling, time slices or snapshots, system states, and representative periods (e.g., days or weeks); comprehensive overviews are given in~\cite{teichgraeber_time-series_2022} and~\cite{hoffmann_review_2020}. Downsampling, one of the simpler methods, lowers the temporal resolution and averages input data over multiple time steps. However, uniformly coarsening the time series can systematically distort investment and operational decisions once variable renewables constitute a significant share of the system~\cite{auer_uncovering_2026}. More refined variants assign different resolutions to individual variables and constraints rather than coarsening uniformly~\cite{gao_fully_2025}. Using representative periods (RPs), on the other hand, is a more complex approach, preserving the original temporal resolution and chronology of the data. RPs are widely used in ESOMs for large-scale energy system studies, such as GenX~\cite{bonaldo_genx_2026}, Tulipa~\cite{soares_siqueira_tulipa_2026}, and LEGO~\cite{wogrin_lego_2022}.

The complexity reduction is achieved by separating the full, chronological, interconnected mathematical problem into distinct periods (e.g., individual 24-hour periods) and solving a reduced number of such periods (e.g., only 7 days of a year). If there were no temporal dynamics in the model, all inaccuracies of the aggregated model would only be caused by the inaccuracy of representing the original days by a reduced number of RPs. But in power and energy systems, such temporal dynamics between individual periods exist, ranging from short-term (e.g., unit-commitment and ramping for thermal generators, intra-day storage with a daily usage pattern, etc.) to long-term (e.g., inter-day storage like pumped hydro, etc.). In RP frameworks, any dynamic exceeding one RP is lost by splitting the full year, unless the relevant constraints are adjusted accordingly.

Therefore, there are two ways to improve solution quality within the representative period framework: Either the selection of RPs is improved (e.g., with better clustering techniques, by choosing more RPs, by including extreme periods, etc.), or the problem is reformulated to capture temporal dynamics better. Obtaining \emph{good} RPs is challenging and part of ongoing research (e.g., see~\cite{wogrin_time_2023, poncelet_selecting_2017, pinto_evaluation_2020}). In this paper, however, we focus on the reformulation of constraints, analyzing the state of the art and identifying shortcomings regarding the original data representation. Most existing work targets improving \textit{inter}-RP linking~\cite{hoffmann_review_2020}, which is concerned with long-term effects spanning multiple RPs. Such approaches typically reintroduce chronology between RPs by adding dedicated inter-period variables and constraints, e.g., to represent seasonal storage~\cite{kotzurTimeSeriesAggregation2018, tejada-arango_enhanced_2018, gonzato_long_2021, pineda_chronological_2018} or to link the operation of short-term storage across periods~\cite{tejada-arangoOpportunityCostIncluding2019}. The handling of \textit{intra}-RP edges, in contrast, has received little attention: existing models treat the period boundaries by either not enforcing the time-linking constraints, fixing the boundary values, or cyclically connecting each RP to itself (see Section~\ref{sec:EdgeHandling}), but none exploit the transition probabilities between RPs. Even recent work that explicitly targets intra-day accuracy does not address this edge problem. Moradi-Sepahvand and Tindemans~\cite{moradi-sepahvand_representative_2024}, for instance, sharpen the representation \emph{within} each representative day through piecewise-linear transitions and an adaptive allocation of more time points to more variable days, but reintroduce inter-day chronology through a dedicated day-mapping construction of the inter-RP type discussed above. The period \emph{edges} themselves are thus still handled no better than by a simple cyclic connection. Our approach instead exploits the transition probabilities between RPs to reformulate the existing boundary constraints directly, without any dedicated extraction or mapping construction.

The main contribution of this paper is a novel reformulation of \textit{intra}-RP constraints, which we call the \Markov. While Transition Matrices have previously been used to link RPs~\cite{tejada-arango_enhanced_2018, wogrin_new_2014, orgaz_temporal_2022}, we instead use the Transition Matrix to supply \textit{expected} boundary values for the existing \textit{intra}-RP constraints, without introducing any additional variables or constraints. The general concept applies to any time-linking constraint, and we also extend it to constraints spanning multiple time steps and to binary variables. Unlike approaches that introduce new transition constructions between periods~\cite{moradi-sepahvand_representative_2024}, our reformulation leaves the number of variables and constraints unchanged. It outperforms existing methods when transitions between periods carry information that current methods discard, while matching them otherwise. The approach is demonstrated using an illustrative case study and the updated IEEE Reliability Test System (RTS-GMLC)~\cite{barrowsIEEEReliabilityTest2020}.

The remainder of this paper is organized as follows. In Section~\ref{sec:EdgeHandling}, the challenge regarding transitions between RPs is introduced. Current state-of-the-art formulations are shown, and their limitations are discussed, after which the suggested solution (the \Markov) is introduced in Section~\ref{sec:Markov}. The approach is extended to connecting multiple time steps (Subsection~\ref{subsec:MultipleTimeSteps}), and a method to handle binary variables is shown (Subsection~\ref{subsec:BinaryVariables}). Lastly, the practical benefits are demonstrated in two case studies and discussed in Section~\ref{sec:CaseStudies}, before concluding and outlining potential future research in Section~\ref{sec:Conclusion}.

\section{Transitions between Representative Periods}\label{sec:EdgeHandling}
When constructing an aggregated model, all \textit{time-linking} constraints (which include variables from at least two different time steps) must be reformulated for the \textit{edges} or \textit{transitions} between those periods. Ensuring this continuity of model states across the time steps that bridge representative periods is a known difficulty when working with representative periods~\cite{pfenningerDealingMultipleDecades2017}. We illustrate three commonly used methods for handling transitions between RPs, using the ramp-up constraint of a thermal power plant (adjusted from~\cite{morales-espana_tight_2013}) as an example. See~\eqref{eq:rampingStandard} for the formulation for consecutive time steps $t\in[1,T]$ and~\eqref{eq:rampingNoEnforcement} for RPs $\rp\in[1,\RP]$ and time steps $k\in[1,K]$ within each $\rp$. In both formulations, the first time step is omitted since there is no \textit{previous} time step for the first one. The variable $p_{t,g}$ denotes the power output at $t$ (or $\rp$ and $k$) of power generator~$g\in\mathcal{G}$, and $\RU_g$ denotes the ramp-up rate of thermal unit $g$ per time step:
\begin{align}
    p_{t,g}-p_{t-1,g} &\leq \RU_g  \nonumber \\
    &\forall t \in \mathcal{T} \setminus \{t_1\},\;g\in \mathcal{G} \label{eq:rampingStandard} \\
    p_{\rp,k,g}-p_{\rp,k-1,g} &\leq \RU_g \nonumber \\
    &\forall \rp\in\mathcal{RP},\; k\in\mathcal{K}\setminus\{k_1\},\;g\in \mathcal{G} \label{eq:rampingNoEnforcement}
\end{align}

\indent \textbf{No Enforcement:} The most straightforward implementation is not enforcing the constraints at the edges of RPs at all. For the Ramp-Up example, this would mean keeping the formulation as shown in~\eqref{eq:rampingNoEnforcement}, where the model can \textit{choose} its variable values for $p_{\rp,k_1,g}$ of each $\rp$ without conforming to any Ramp-Up restrictions. This can also lead to significant errors for short-term storage, as the model can assume full storage at the beginning of each RP and deplete it during the RP, essentially leading to free energy. That being said, if $K$ is sufficiently large, the introduced error is expected to be low (e.g., for an hourly model of a year where $T=K=8760$).

\textbf{Fixed Initial Values:} To have a more conservative approach, initial values can be provided for those time steps at the edges of RPs. For the Ramp-Up example, one can fix $p_{\rp,k_1,g}$ using some pre-determined fixed output $\FP_{\rp,g}$, which the modeler has to provide. That means one needs $\FP_{\rp,g}\;\forall g\in\mathcal{G},\rp\in\mathcal{RP}$. In most cases, this is not practical, since the optimal boundary values are not known a priori. Example models using this method are listed in Table~\ref{tab:edgeMethods}: \cite{wogrin_lego_2022}, \cite{gonzato_long_2021}, \cite{brown_pypsa_2018}, \cite{schwele_unit_2020}, \cite{soares_siqueira_tulipa_2026}, \cite{pfenninger_calliope_2018}.

\textbf{Cyclic Connection:} Finally, the constraints can be constructed to connect the edges of each $\rp$ with itself cyclically. For the Ramp-Up example, that would mean using $k_K$ as the predecessor of $k_1$, i.e., $p_{\rp,k_1,g}-p_{\rp,k_K,g}\leq\RU_g\;\forall g\in\mathcal{G},\rp\in\mathcal{RP}$, therefore also enforcing a constraint for $k_1$. This reduces the freedom of the model without requiring a priori knowledge about system states, which makes it a popular choice in state-of-the-art models, e.g.: \cite{bonaldo_genx_2026}, \cite{wogrin_lego_2022}, \cite{kotzurTimeSeriesAggregation2018}, \cite{brown_pypsa_2018}, \cite{soares_siqueira_tulipa_2026}, \cite{gabrielli_optimal_2018}, \cite{pfenninger_calliope_2018}, \cite{helisto_backbone_2019}, \cite{johnston_switch_2019} (see Table~\ref{tab:edgeMethods}).

While superior to the other two discussed approaches, the \Cyclic assumes that each RP always follows itself. It is easy to show that this can also lead to significant modeling errors: Figure~\ref{fig:TransitionPeriod} illustrates this for a Ramp-Up example. Assume that for the actual optimal solution, $p_{\rp_1, k_1, g}$ would be at a \textit{high} value and $p_{\rp_1, k_K, g}$ at a \textit{low} value. If $\rp_2$ follows $\rp_1$, $p_{\rp_2, k_1, g}$ should optimally also start on a \textit{low} level and continue from there. However, in this method, $p_{\rp_1,k_1,g}$ is bound to $p_{\rp_1,k_K,g}$ and therefore forces $p_{\rp_1,k_K,g}$ to deviate from the optimal result, creating an error since the cyclic Ramp-Up constraint would otherwise be violated.

\begin{table}[!t]
\caption{Edge-handling methods used in state-of-the-art models.}
\label{tab:edgeMethods}
\centering
\resizebox{\columnwidth}{!}{%
\begin{tabular}{@{}lll@{}}
\hline
Model & Edge-handling method & Note \\ \hline
GenX~\cite{bonaldo_genx_2026} & Cyclic Connection & Cyclic intra-period SoC \\
Tulipa~\cite{soares_siqueira_tulipa_2026} & Fixed Initial Values or Cyclic & Non-draining if values are fixed \\
LEGO~\cite{wogrin_lego_2022} & Fixed Initial Values or Cyclic & Fixed if $\RP=1$, else Cyclic \\
Kotzur et al.~\cite{kotzurTimeSeriesAggregation2018} & Cyclic Connection & For each time frame \\
Gonzato et al.~\cite{gonzato_long_2021} & Fixed Initial Values & Via inter-period constraints \\
PyPSA~\cite{brown_pypsa_2018} & Fixed Initial Values or Cyclic & $\RP=1$, not a classical RP framework \\
Schwele et al.~\cite{schwele_unit_2020} & Fixed Initial Values & Fixed initial output/status per RP \\
Calliope~\cite{pfenninger_calliope_2018} & Fixed Initial Values or Cyclic & Cyclic is the default \\
Gabrielli et al.~\cite{gabrielli_optimal_2018} & Cyclic Connection & Periodicity constraints \\
Backbone~\cite{helisto_backbone_2019} & Cyclic Connection & Cyclic node-state boundaries \\
Switch~\cite{johnston_switch_2019} & Cyclic Connection & Circular within each RP \\ \hline
\end{tabular}}
\end{table}

\begin{figure}
\centering
\scalebox{0.7}{
\begin{tikzpicture} 
\begin{axis} [
    scale=0.95,
    legend style={font=\small, cells={align=left}},
    legend cell align=left,
    legend pos=south east,
    ymin=0,ymax=1,
    ylabel={$p_{\rp,k,g_1}$},
    xtick={1, 12, 24, 25, 36, 48},
    xticklabels={$k$=1, 12, 24 \kern1mm, \kern1mm 1, 12, 24},
    axis lines=left,
    height=5cm,width=\textwidth/2,
    clip=false,]
    \addplot [
        no marks, blue
    ]
    coordinates{ 
        (1, 1) 
        (5, 0.55) 
        (8, 0.6) 
        (15, 0.5) 
        (17, 0.45) 
        (20, 0.2) 
        (24, 0.0) 
        (26, 0.1) 
        (31, 0.2) 
        (39, 0.3) 
        (41, 0.4) 
        (44, 0.3) 
        (48, 0.4) 
    };
    \addplot [
        no marks, red, very thick, dotted
    ]
    coordinates{ 
        (25, 1) 
        (29, 0.55) 
        (32, 0.6) 
        (39, 0.5) 
        (41, 0.45) 
        (44, 0.2) 
        (48, 0.0) 
    };
    \node[red, align=left] (text) at (axis cs:34,0.8){Constraint\\violation};
    \node (source) at (axis cs:24,0){};
    \node (destination) at (axis cs:25,1){};
    \draw[red, line width=2pt, arrows={-Triangle[angle=80:8pt,red,fill=red]}] (source) -- (destination);

    \draw[{| Triangle[]-Triangle[] |},line width=1pt] (axis cs:1,-0.23) -- node[fill=white,sloped]{$\rp_1$}(axis cs:24,-0.23);
    \draw[{| Triangle[]-Triangle[] |},line width=1pt] (axis cs:25,-0.23) -- node[fill=white,sloped]{$\rp_2$}(axis cs:48,-0.23);
    \legend{Optimal $p$,Cyclic $p_{\rp_1}$}
\end{axis} 
\end{tikzpicture}
}
\caption{Illustration of a Transition Period $\rp_1$ where the theoretically optimal $p$ cannot be obtained due to a \Cyclic of the Ramp-Up constraint.}
\label{fig:TransitionPeriod}
\end{figure}

\section{The Markov Transition Method}\label{sec:Markov}
This section introduces an innovative approach, \Markov, leveraging the Transition Matrix to overcome the shortcomings of traditional approaches by reintroducing information about the transitions between periods. The Transition Matrix $P$ is constructed from the original chronological data by counting the number of transitions $N_{\rp',\rp}$ from $\rp'$ to $\rp$. Since the reformulated constraints introduced below require the \textit{expected previous} time step of each $\rp$, the counts are normalized over the transitions entering $\rp$, i.e., $P_{\rp',\rp}=N_{\rp',\rp}/\sum_{\rp''\in\mathcal{RP}}N_{\rp'',\rp}$ is the probability that a given $\rp$ was preceded by $\rp'$ (the transition matrix of the time-reversed sequence). A visual example of such a Transition Matrix is shown in Fig.~\ref{fig:TransitionMatrix}. The idea of deriving transition probabilities between clustered states is not new: it expands the \emph{system states} approach introduced over a decade ago~\cite{wogrin_new_2014} and still applied in recent temporal-aggregation work~\cite{orgaz_temporal_2022}. There, however, the transitions serve to \emph{construct the aggregated model itself}, adding dedicated state variables and constraints. In contrast, we use the Transition Matrix only to supply \emph{expected} boundary values for the existing \textit{intra}-RP constraints, leaving the number of variables and constraints unchanged. None of the state-of-the-art models surveyed in Table~\ref{tab:edgeMethods} exploit this transition information at their period boundaries; they all rely on the traditional edge-handling approaches of Section~\ref{sec:EdgeHandling}.

\begin{figure}
\centering
\begin{tikzpicture}[shorten >=1pt, auto,
    node distance=2.5cm, scale=0.7, 
    transform shape,
    block/.style={draw,rectangle,minimum width={width("$t_{145}-t_{168}$")+10pt}}]
    
    \node(Title) at (0,6.5) {Original Time Series};
    \node[block, fill=color1, text=white] (A) at (0,6) {$t_{1}-t_{24}$};
    \node[block, fill=color2] (B) at (0,5) {$t_{25}-t_{48}$};
    \node[block, fill=color1, text=white] (C) at (0,4) {$t_{49}-t_{72}$};
    \node[block, fill=color1, text=white] (D) at (0,3) {$t_{73}-t_{96}$};
    \node[block, fill=color2] (E) at (0,2) {$t_{97}-t_{120}$};
    \node[block, fill=color3, text=white] (F) at (0,1) {$t_{121}-t_{144}$};
    \node[block, fill=color3, text=white] (G) at (0,0) {$t_{145}-t_{168}$};
    
    \path [->] (A) edge node[left] {} (B);
    \path [->] (B) edge node[left] {} (C);
    \path [->] (C) edge node[left] {} (D);
    \path [->] (D) edge node[left] {} (E);
    \path [->] (E) edge node[left] {} (F);
    \path [->] (F) edge node[left] {} (G);
    \path [->] (G) edge[bend left=90] node[left] {} (A);

    \node[block, fill=color1, text=white, minimum width=30pt] (RP1) at (3,6) {$t_x-t_y$};
    \node[below right=-0.5cm and 0cm of RP1]{Time steps assigned to $\rp_1$};
    \node[block, fill=color2, minimum width=30pt, below=0.1cm of RP1] (RP2) {$t_x-t_y$};
    \node[below right=-0.5cm and 0cm of RP2]{Time steps assigned to $\rp_2$};
    \node[block, fill=color3, text=white, minimum width=30pt, below=0.1cm of RP2] (RP3) {$t_x-t_y$};
    \node[below right=-0.5cm and 0cm of RP3]{Time steps assigned to $\rp_3$};

    \node [shape=rectangle, draw, align=left, below right=1.2cm and -1.3cm of RP3](table) {
        Transition Matrix \\
 \begin{tabular}{lllll}
      &                            &       \multicolumn{3}{c}{To $\rp$:}       \\
      &                            & \cellcolor{color1}\color{white}$\rp_1$ & \cellcolor{color2}$\rp_2$ & \cellcolor{color3}\color{white}$\rp_3$ \\ \cline{3-5}
      & \multicolumn{1}{l|}{\cellcolor{color1}\color{white}$\rp_1$} & $1/3$   & $1$   &       \\
From $\rp'$: & \multicolumn{1}{l|}{\cellcolor{color2}$\rp_2$} & $1/3$   &       & $1/2$   \\
      & \multicolumn{1}{l|}{\cellcolor{color3}\color{white}$\rp_3$} & $1/3$   &       & $1/2$
\end{tabular}
    };
\end{tikzpicture}
\caption{Illustration of the Transition Matrix for a 168h period, which is represented by 3 RPs with 24 hours each. Each cell gives the probability $P_{\rp',\rp}$ that the RP in the column was preceded by the RP in the row, so each column sums to 1.}
\label{fig:TransitionMatrix}
\end{figure}

\subsection{Adjusting Constraints with Transition Probabilities}
Time-linking constraints are altered to enhance the \Cyclic by replacing the values of the \textit{previous} time step of the \textit{same} $\rp$ with the values of the \textit{expected previous} time step. Note that there are methods to reconstruct some chronology for seasonal storage (e.g., \cite{tejada-arango_enhanced_2018}, \cite{gonzato_long_2021}) by creating additional variables and constraints to improve \textit{inter}-period links. However, our method aims to improve current \textit{intra}-RP formulations without increasing the number of variables or constraints. 

See~\eqref{eq:rampingMarkov} for an example of calculating the value of the \textit{expected previous} time step for the Ramp-Up constraint, where an \textit{expected value} of the variable in question is calculated. This also works similarly when calculating \textit{expected next} time steps, then normalizing the counts $N$ over the transitions leaving each period instead.
\begin{multline}
    p_{\rp,k_0,g}\coloneqq \sum_{\rp'}(p_{\rp',k_K,g}\cdot P_{\rp', \rp})\\
    \forall \rp\in\mathcal{RP},\;g\in\mathcal{G}\label{eq:rampingMarkov}
\end{multline}
Here, $k_0$ denotes the predecessor of the first time step $k_1$ of an $\rp$: for the \Cyclic it equals $k_K$ of the same $\rp$, whereas for the \Markov it is the \textit{expected previous} value of~\eqref{eq:rampingMarkov}. The Ramp-Up restriction at the period boundary is then enforced as
\begin{equation}
    p_{\rp,k_1,g}-p_{\rp,k_0,g}\leq \RU_g \quad \forall \rp\in\mathcal{RP},\;g\in\mathcal{G}.\label{eq:rampingMarkovEdge}
\end{equation}
Using an \textit{expected} boundary value is achieved by replacing the unknown predecessor value by the probability-weighted average $\sum_{\rp'}(p_{\rp',k_K,g}\cdot P_{\rp',\rp})$. This enforces the constraint \textit{in expectation} over the possible predecessors of each $\rp$, while leaving the number of variables and constraints unchanged.

Similarly, \textit{intra}-RP storage formulations can also be improved by the Transition Matrix, where a unit $s$ with storage level $l_{\rp, k, s}$ can generate (i.e., discharge) with $p_{\rp, k, s}$ and consume (i.e., charge) with $c_{\rp, k, s}$. Constraint~\eqref{eq:storage:a} shows the standard intra-RP constraint for the storage level update, and~\eqref{eq:storage:b} shows the special case for the first time step, which relies on the described \textit{expected} value\footnote{Note that $\Delta_k$, the duration of one time step, is used to convert $c$ and $p$ (originally in unit of power) to energy. Efficiencies, which are commonly used in storage constraints, are left out for brevity.}:
\begin{subequations}\label{eq:storage}
    \begin{multline}
    l_{\rp,k,s} =  l_{\rp,k-1,s} + (c_{\rp,k,s}- p_{\rp,k,s})\cdot\Delta_k\\
    \hspace{2em}\forall \rp\in\mathcal{RP},\; k\in\mathcal{K}\setminus\{k_1\},\; s\in\mathcal{S}\label{eq:storage:a}
    \end{multline}
    \begin{multline}
    l_{\rp,k_1,s} = \sum_{\rp'}(l_{\rp',k_K,s}\cdot P_{\rp', \rp})+ (c_{\rp,k_1,s} - p_{\rp,k_1,s})\cdot\Delta_k \\
    \forall \rp\in\mathcal{RP},\;s\in\mathcal{S}\label{eq:storage:b}
    \end{multline}
\end{subequations}

\subsection{Connecting Multiple Time Steps}\label{subsec:MultipleTimeSteps}
The approach can also handle constraints that connect multiple time steps. We demonstrate it using a Minimum Up-Time constraint, adjusted from~\cite{morales-espana_tight_2013}. Constraints~\eqref{eq:MinUpTime:a} and~\eqref{eq:MinUpTime:b} govern the relationship between the production $p_{\rp, k, g}$, the maximum production capability $\PC_g$, and the binary variables $u_{\rp, k, g}$ (unit is committed), $\su_{\rp, k, g}$ (unit is starting up), and $\sd_{\rp, k, g}$ (unit is shutting down). Constraint~\eqref{eq:MinUpTime:c} ensures that a unit stays committed for at least $\MU_g$ time steps after a start-up.

The key difference to the single-step case (e.g., Ramp-Up) is that the sliding window $\sum_{k'=k-\MU_g+1}^{k}$ in~\eqref{eq:MinUpTime:c} reaches back further than one time step. For the first time steps of an $\rp$ (i.e., for $k_1,\dots,k_{\MU_g-1}$), the lower summation index $k-\MU_g+1$ becomes non-positive and thus refers to time steps \textit{preceding} the period, which we denote by indices $k\leq 0$ (with $k_0$ being the immediate predecessor of $k_1$, $k_{-1}$ the one before, and so on). Under the \Cyclic, these pre-period values are taken from the same $\rp$ by wrapping the index around, i.e., replacing any $k\leq 0$ by $k+K$. Under the \Markov, each such pre-period value is instead replaced by its \textit{expected value} over the possible predecessors $\rp'$, weighted by the transition probabilities $P_{\rp',\rp}$, analogously to~\eqref{eq:rampingMarkov}. Definition~\eqref{eq:MinUpTimeMarkov} gives these expected values for $\su_{\rp, k, g}$ for all non-positive offsets $k\leq 0$.
\begin{subequations}\label{eq:MinUpTime}
    \begin{multline}
        p_{\rp,k,g}\leq \PC_g\cdot u_{\rp, k, g}\\
        \forall \rp\in\mathcal{RP},\; k\in\mathcal{K},\; g\in\mathcal{G}\label{eq:MinUpTime:a}
    \end{multline}
    \begin{multline}
        u_{\rp, k, g} - u_{\rp, k-1, g} = \su_{\rp, k, g} - \sd_{\rp, k, g}\\
        \forall \rp\in\mathcal{RP},\; k\in\mathcal{K},\; g\in\mathcal{G}\label{eq:MinUpTime:b}
    \end{multline}
    \begin{multline}
        \sum_{k'=k-\MU_g + 1}^{k}\su_{\rp, k', g} \leq u_{\rp, k, g}\\
        \forall \rp\in\mathcal{RP},\; k\in\mathcal{K},\; g\in\mathcal{G}\label{eq:MinUpTime:c}
    \end{multline}
\end{subequations}
\begin{multline}
    \su_{\rp, k, g} \coloneqq \sum_{\rp'}(\su_{\rp', k+K, g}\cdot P_{\rp', \rp}) \\
    \forall \rp\in\mathcal{RP},\; k\leq 0,\; g\in\mathcal{G}\label{eq:MinUpTimeMarkov}
\end{multline}

\subsection{Dealing with Binary Variables}\label{subsec:BinaryVariables}
When calculating \textit{expected values} for binary variables using our proposed approach (e.g., with~\eqref{eq:MinUpTimeMarkov}), one may obtain non-binary results. Assume transition probabilities $P_{\rp_1, \rp_1}=0.7$, $P_{\rp_2,\rp_1}=0.3$ and unit commitments $u_{\rp_1,k_K,g}=0$ and $u_{\rp_2,k_K,g}=1$. The \textit{expected previous} of $u_{\rp_1,k_1,g}$ would then be $u_{\rp_1,k_0,g}=0.3$, which is fractional (i.e., $\notin\{0,1\}$). However, $u_{\rp_1,k_1,g}$ can only differ from $u_{\rp_1,k_0,g}$ by $\su_{\rp_1,k_1,g}-\sd_{\rp_1,k_1,g}\in\{-1,0,1\}$ due to~\eqref{eq:MinUpTime:b}, since all other variables in that equation are binary. Starting from a fractional $u_{\rp_1,k_0,g}$, no integer value of $u_{\rp_1,k_1,g}$ can satisfy~\eqref{eq:MinUpTime:b}, so the fractional value would lead to infeasibility.

To leverage the \Markov approach while also enforcing the binary nature of such variables as much as possible, we relax\footnote{\emph{Relaxing} in the sense of allowing decimal values between $0$ and $1$.} all binaries at time steps which are directly influenced by \textit{expected values}. This can range from single time steps (e.g., Ramp-Up) to multiple time steps (e.g., Minimum Up-Time). The time span for the relaxation of an arbitrary binary variable $v_{\rp, k, *}$ (with indices $\rp,k$ and potentially additional indices '$*$') depends on its indices. For example, the Minimum Up-Time can be different for generators $g_1$ and $g_2$, and the relaxation must be adjusted accordingly. Importantly, only the variables \textit{within} the relevant Minimum Up-/Down-Time window are relaxed. All remaining time steps keep their binary domain. The relaxed variables \textit{can}, but do not have to, take non-binary values, which would require a correction step (e.g., if the solution should be used as a dispatch schedule). In practice, this affects only a small share of variables, as we observe in our numerical experiments in the next section.

\section{Case Studies and Discussion}\label{sec:CaseStudies}
\subsection{Experimental Setup}\label{subsec:ExperimentalSetup}
We present two case studies: The first, which we solve as an operational problem, is specifically designed to show the benefits and limitations of the \Markov, and is described in Subsection~\ref{subsec:illustrativeCaseStudy}. The second, larger, and more realistic case is based on the updated IEEE Reliability Test System (RTS-GMLC)~\cite{barrowsIEEEReliabilityTest2020}, described in Subsection~\ref{subsec:realisticCaseStudy}. 

For both, we first cluster the full chronological data by a given number of RPs. We create a \Reference by reconstructing the full modeling horizon from this clustered data, where we use copies of the representatives that have been identified to replace each of the original days. A visualization of this for three representative days can be seen in Fig.~\ref{fig:TruthCase}. The aim is to remove the aggregation error from our comparisons later. Otherwise, if we compared the results of our edge-handling candidates to the results of the \textit{original, unclustered data}, the comparison would also include an error that stems from an inaccurate representation of the data. With our \Reference there is no aggregation error, since if we ran the clustering again, the exact same clusters would be found and each day would be represented perfectly by its RP (since it is identical). 

We have applied our framework to constraints for Ramping (constraints~\crefrange{eq:rampingNoEnforcement}{eq:rampingMarkovEdge} for Ramp-Up, similarly also for Ramp-Down), Minimum Up-/Down-Time (constraints~\eqref{eq:MinUpTime} and definition~\eqref{eq:MinUpTimeMarkov} for Minimum Up-Time, similarly also for Minimum Down-Time), and Intra-Day Storage (constraints~\eqref{eq:storage}). 

All methods are implemented in the open-source LEGO model~\cite{wogrin_lego_2022} and solved with \emph{Gurobi 13.0.0} using a relative MIP gap of $0.3\%$\footnote{To handle the immense RAM requirements of the \Reference models, \textit{NodefileStart} is set to 100~GB to offload some memory from RAM (if possible), and threads are limited to 10 per Gurobi instance.}. A Windows server with an Intel Xeon Gold processor at 3.6~GHz and 2~TB of RAM is used to solve the models, and we report computational effort in \textit{work units}\footnote{Work units are a solver-internal, deterministic measure of computational effort, independent of machine load and other external factors.}. Code and data are available at \url{https://github.com/IEE-TUGraz/LEGO-Pyomo/tree/research/MarkovTransition}.

\begin{figure}
\centering
\begin{tikzpicture}[shorten >=1pt, auto,
    node distance=2.5cm, scale=0.55, 
    transform shape,
    block/.style={draw,rectangle,minimum width={width("$t_{145}-t_{168}$")+10pt}, minimum height=14pt}]
    
    \node[block,draw=white] (T) at (0,6.7) {Orig. Chronological Model};
    \node[block, fill=color1, text=black, fill opacity=0.4, text opacity=1.0] (A) at (0,6) {$t_{1}-t_{24}$};
    \node[block, fill=color2, fill opacity=0.3, text opacity=1.0] (B) at (0,5) {$t_{25}-t_{48}$};
    \node[block, fill=color1, text=white] (C) at (0,4) {$t_{49}-t_{72}$};
    \node[block, fill=color1, text=black, fill opacity=0.4, text opacity=1.0] (D) at (0,3) {$t_{73}-t_{96}$};
    \node[block, fill=color2] (E) at (0,2) {$t_{97}-t_{120}$};
    \node[block, fill=color3, text=white] (F) at (0,1) {$t_{121}-t_{144}$};
    \node[block] (G) at (0,0) {...};
    
    \path [->] (A) edge node[left] {} (B);
    \path [->] (B) edge node[left] {} (C);
    \path [->] (C) edge node[left] {} (D);
    \path [->] (D) edge node[left] {} (E);
    \path [->] (E) edge node[left] {} (F);
    \path [->] (F) edge node[left] {} (G);
    
    \node[block,draw=white] (T) at (4,6.7) {Reference Case};
    \node[block, fill=color1, text=white] (A1) at (4,6) {$t_{49}-t_{72}$};
    \node[block, fill=color2] (B1) at (4,5) {$t_{97}-t_{120}$};
    \node[block, fill=color1, text=white] (C1) at (4,4) {$t_{49}-t_{72}$};
    \node[block, fill=color1, text=white] (D1) at (4,3) {$t_{49}-t_{72}$};
    \node[block, fill=color2] (E1) at (4,2) {$t_{97}-t_{120}$};
    \node[block, fill=color3, text=white] (F1) at (4,1) {$t_{121}-t_{144}$};
    \node[block] (G1) at (4,0) {...};

    \path [->] (A1) edge node[left] {} (B1);
    \path [->] (B1) edge node[left] {} (C1);
    \path [->] (C1) edge node[left] {} (D1);
    \path [->] (D1) edge node[left] {} (E1);
    \path [->] (E1) edge node[left] {} (F1);
    \path [->] (F1) edge node[left] {} (G1);

    \contourlength{0.4pt}
    \node[block, draw=white, minimum width=30pt] (P1) at (-3,6) {\textcolor{gray}{Represented by} \contour{black}{\textcolor{color1}{\boldmath$\rp_1$}}};
    \path [->] (P1) edge[gray] node[left] {} (A);

    \contourlength{0.6pt}
    \node[block, draw=white, minimum width=30pt] (P2) at (-3,5) {\textcolor{gray}{Represented by} \contour{black}{\textcolor{color2}{\boldmath$\rp_2$}}};
    \path [->] (P2) edge[gray] node[left] {} (B);

    \contourlength{0.4pt}
    \node[block, draw=white, minimum width=30pt] (RP1) at (-2.8,4) {Identified as \contour{black}{\textcolor{color1}{\boldmath$\rp_1$}}};
    \path [->] (RP1) edge node[left] {} (C);

    \node[block, draw=white, minimum width=30pt] (P4) at (-3,3) {\textcolor{gray}{Represented by} \contour{black}{\textcolor{color1}{\boldmath$\rp_1$}}};
    \path [->] (P4) edge[gray] node[left] {} (D);

    \contourlength{0.6pt}
    \node[block, draw=white, minimum width=30pt] (RP2) at (-2.8,2) {Identified as \contour{black}{\textcolor{color2}{\boldmath$\rp_2$}}};
    \path [->] (RP2) edge node[left] {} (E);

    \contourlength{0.4pt}
    \node[block, draw=white, minimum width=30pt] (RP3) at (-2.8,1) {Identified as \contour{black}{\textcolor{color3}{\boldmath$\rp_3$}}};
    \path [->] (RP3) edge node[left] {} (F);
    
    \draw [->, color1, thick] (C.east) to[out=0, in=180] (A1.west);
    \draw [->, color2, thick] (E.east) to[out=0, in=180] (B1.west);
    \path [->, color1, thick] (C) edge node[below] {} (C1);
    \draw [->, color1, thick] (C.east) to[out=0, in=180] (D1.west);
    \path [->, color2, thick] (E) edge node[near start] {} (E1);
    \path [->, color3, thick] (F) edge node[left] {} (F1);

    \node[block,draw=white] (T) at (7,6.7) {Cyclic Connection};
    \node[block, fill=color1, text=white] (RP1C) at (7,6) {$t_{49}-t_{72}$};
    \node[block, fill=color2] (RP2C) at (7,5) {$t_{97}-t_{120}$};
    \node[block, fill=color3, text=white] (RP3C) at (7,4) {$t_{121}-t_{144}$};
    
    \draw [->] (RP1C.east)to[out=315, in=225, looseness=1.5] (RP1C.west);
    \draw [->] (RP2C.east)to[out=315, in=225, looseness=1.5] (RP2C.west);
    \draw [->] (RP3C.east)to[out=315, in=225, looseness=1.5] (RP3C.west);

    \node[block,draw=white] (T) at (10,6.7) {No Enforcement};
    \node[block, fill=color1, text=white] (RP1N) at (10,6) {$t_{49}-t_{72}$};
    \node[block, fill=color2] (RP2N) at (10,5) {$t_{97}-t_{120}$};
    \node[block, fill=color3, text=white] (RP3N) at (10,4) {$t_{121}-t_{144}$};
    
    \draw [|-] (RP1N.west)+(-4mm,0) -- +(0.5mm,0);
    \draw [|-] (RP2N.west)+(-4mm,0) -- +(0.5mm,0);
    \draw [|-] (RP3N.west)+(-4mm,0) -- +(0.5mm,0);
    \draw [-|] (RP1N.east) -- +(4mm,0);
    \draw [-|] (RP2N.east) -- +(4mm,0);
    \draw [-|] (RP3N.east) -- +(4mm,0);

    \node[block,draw=white] (T) at (8.5,3) {Markov Transition};
    \node[block, fill=color1, text=white] (RP1M) at (8.5,1.7) {$t_{49}-t_{72}$};
    \node[block, fill=color2] (RP2M) at (7,0.3) {$t_{97}-t_{120}$};
    \node[block, fill=color3, text=white] (RP3M) at (10,0.3) {$t_{121}-t_{144}$};
    
    \draw [->] (RP1M.east) to[out=45, in=135, looseness=1.5] node [above] {\small $P_{\rp_1,\rp_1}$} (RP1M.west);
    \draw [->] (RP1M) to[bend right] node [left] {\small $P_{\rp_1,\rp_2}$} (RP2M);
    \draw [->] (RP1M) to[bend left] node [right] {\small $P_{\rp_1,\rp_3}$} (RP3M);
    \draw [->, dashdotted] (RP2M.south east) to[out=315, in=225, looseness=1.5] (RP2M.south west);
    \draw [->, dashdotted] (RP2M.north)+(2mm,0) to[bend right] (RP1M.south west);
    \draw [->, dashdotted] (RP2M.east) to[bend right, looseness=1.5] (RP3M.west);
    \draw [->, dashdotted] (RP3M.south east) to[out=315, in=225, looseness=1.5] (RP3M.south west);
    \draw [->, dashdotted] (RP3M.north)+(-2mm,0) to[bend left] (RP1M.south east);
    \draw [->, dashdotted] (RP3M.west) to[bend right, looseness=1.5] (RP2M.east);

\end{tikzpicture}
\caption{Illustration of the creation of the \Reference for 3 RPs, where time steps of the identified RPs replace the original time steps. In contrast, the \Markov, \Cyclic, and \NoEnf only use each RP once.}
\label{fig:TruthCase}
\end{figure}

\subsection{Illustrative Case Study}\label{subsec:illustrativeCaseStudy}
\begin{figure}[]
\centering
\includegraphics[width=1.0\columnwidth]{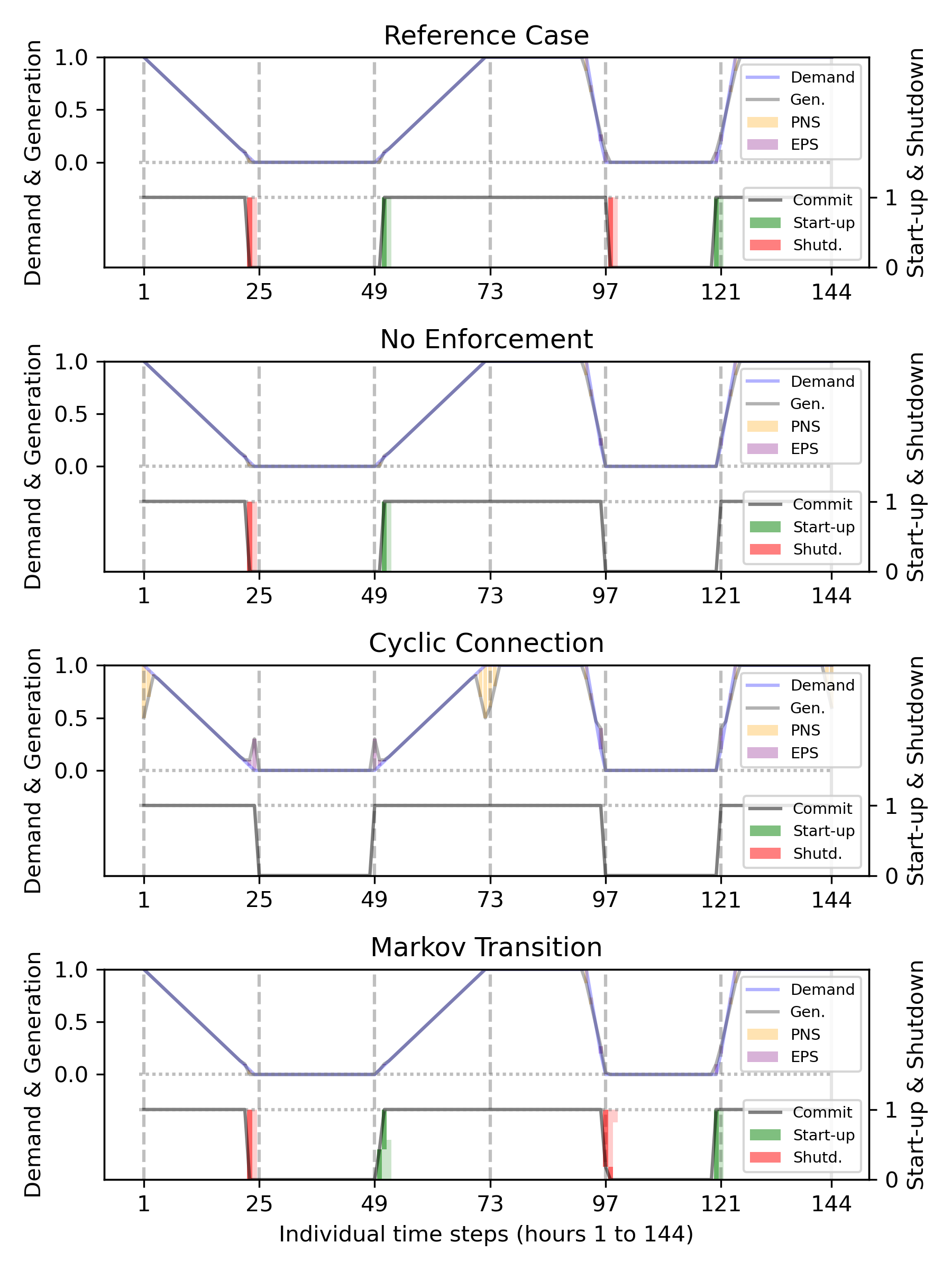}
\caption{Operational decisions for the illustrative case. All power values are normalized.}
\label{fig:illustrative}
\end{figure}

The power system for the illustrative case study consists of a single bus, one thermal generator with unit commitment and ramping, and six consecutive representative days. Each of the six days forms its own RP, so the Transition Matrix contains the exact predecessor of every period, which is the best case for the \Markov. We solve it as an operational problem without investments. The \Markov approximates the objective function value of the \Reference better than \NoEnf and \Cyclic. It only deviates by 8\%, compared to 15\% and 150\%, respectively. More importantly, only the \Markov can exactly determine all start-up and shutdown procedures (see Table~\ref{tab:caseStudy} for exact results). Figure~\ref{fig:illustrative} shows the actual operating decisions for the thermal generator responding to the demand throughout the six days, where the top half shows the normalized demand, generation, power-not-served (PNS), and excess-power-served (EPS). The bottom half shows the binary decisions of commitment, start-up, and shutdown for the generator. One can also see the non-binary values in hours 49--51 and 96--98 for the \Markov due to the relaxation of binaries within the Minimum Up-/Down-Time window (as described in Subsection~\ref{subsec:BinaryVariables}).

\begin{table}[]
\caption{Results of the illustrative case study.}
\label{tab:caseStudy}
\centering
\resizebox{\columnwidth}{!}{%
\begin{tabular}{lrrc}
\hline 
Edge-Handling Method & \makecell{Objective Function\\Value (M EUR)} & Solve Time & \makecell{\# Start-ups /\\ Shutdowns} \\ \hline

Reference Case & $27.71$ & $<0.3$s & $2$ / $2$ \\

No Enforcement & \textcolor{YellowOrange}{$-15\%$} & $<0.3$s & \textcolor{YellowOrange}{$-50\%$} / \textcolor{YellowOrange}{$-50\%$}\\ 
 
Cyclic Connection & \textcolor{BrickRed}{$+150\%$} & $<0.3$s & \textcolor{BrickRed}{$-100\%$} / \textcolor{BrickRed}{$-100\%$}\\ 

Markov Transition & \textcolor{ForestGreen}{$-8\%$} & $<0.3$s & \textcolor{ForestGreen}{$\pm0\%$} / \textcolor{ForestGreen}{$\pm0\%$} \\ \hline

\multicolumn{4}{l}{\footnotesize \makecell[l]{Percentages relative to the reference case. Closest result in \textcolor{ForestGreen}{green}, second in \textcolor{YellowOrange}{orange}\\and worst in \textcolor{BrickRed}{red}. Coloring omitted where values are negligibly small.}}
\end{tabular}}
\end{table}

\subsection{Realistic Case Study}\label{subsec:realisticCaseStudy}
For the realistic case study, we solve the RTS-GMLC test system, which aims to include technologies and configurations that could be encountered in any system~\cite{barrowsIEEEReliabilityTest2020}. It was deliberately chosen so that the full chronological model (serving as the \Reference) can still be solved in reasonable time. The dataset consists of 73 buses in 3 regions, connected by 121 transmission lines (modeled as a Transport Problem), 73 thermal generators, 62 renewable energy sources, and 21 storage units (all modeled as intra-RP storage). We solve it as an operational problem (see Subsection~\ref{subsubsec:OperationalAssessment}) as well as an investment problem (see Subsection~\ref{subsubsec:InvestmentAssessment}). For the RP framework, we solve it for 3, 5, and 7 representative days (in hourly resolution). Additionally, we add variations where we scale the demand's variability: each hourly deviation from the annual mean demand is reduced to 50\%, 70\%, and 90\% of its original value, flattening peaks and valleys symmetrically while keeping the total annual demand unchanged (totaling 4 demand-variability levels, including the original).

\begin{figure}[]
\centering
\includegraphics[width=1\columnwidth]{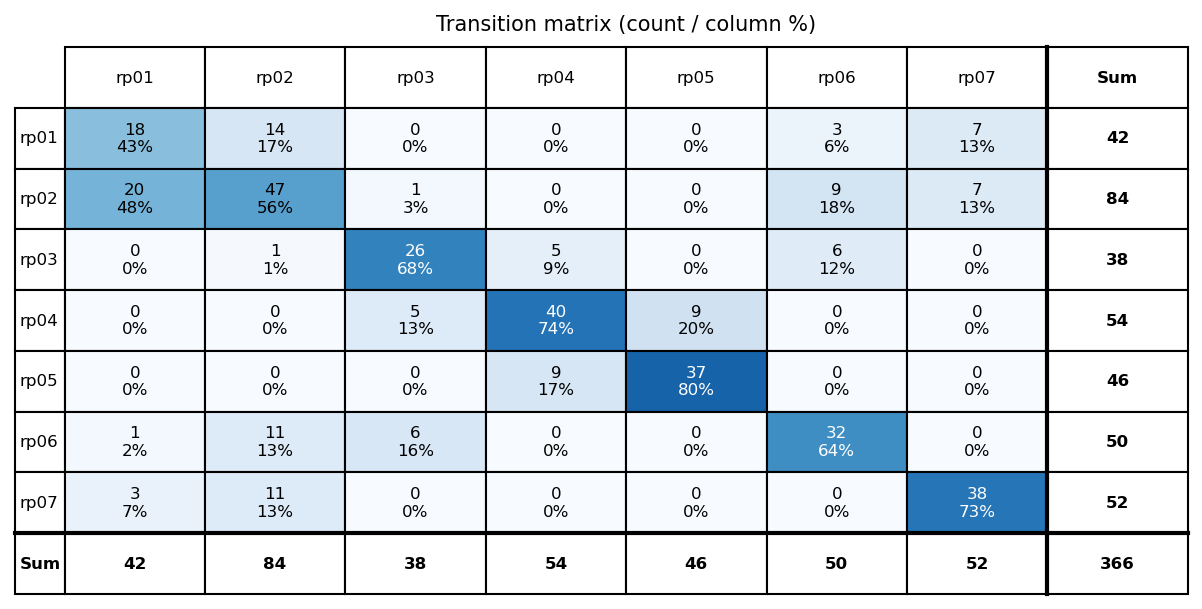}
\caption{Transition counts $N_{\rp',\rp}$ for RTS-GMLC with 7 RPs, from the RP in the rows to the RP in the columns, together with the probabilities $P_{\rp',\rp}$ (each column sums to 100\%, see Section~\ref{sec:Markov}). The background color assists as a visual guide, where darker shades indicate higher values.}
\label{fig:TransitionMatrix7RPs}
\end{figure}

\begin{figure}[]
\centering
\includegraphics[width=1\columnwidth]{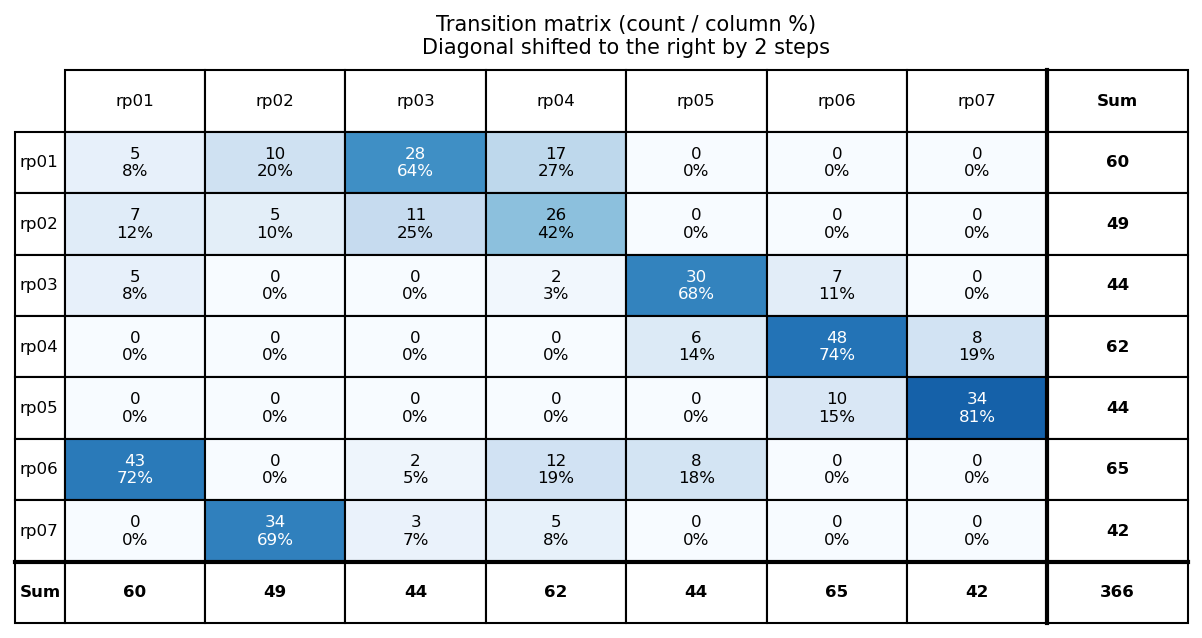}
\caption{Transition counts and probabilities (as in Fig.~\ref{fig:TransitionMatrix7RPs}) for RTS-GMLC with 7 RPs with a diagonal that is shifted to the right by 2 steps.}
\label{fig:TransitionMatrix7RPsShift2}
\end{figure}

The structure of the Transition Matrix has a significant impact on the performance of the different edge-handling methods compared to each other (as will be discussed in Subsection~\ref{subsec:Discussion}). Therefore, we create three variants of the RTS-GMLC dataset to investigate this impact: the original dataset (whose Transition Matrix is close to a diagonal matrix, see Fig.~\ref{fig:TransitionMatrix7RPs}), and two additional variants where we shift the diagonal of the Transition Matrix by 1 and 2 steps, respectively\footnote{For these variants, each row of the transition-count matrix is shifted cyclically, and a new sequence for a year of (unchanged) RPs is sampled from the shifted matrix. The \Reference and all edge-handling methods of a variant are then built from this new sequence.}. See Fig.~\ref{fig:TransitionMatrix7RPsShift2} for an example with a shift of 2 steps.

\subsubsection{Operational Assessment}\label{subsubsec:OperationalAssessment}
\begin{figure}[]
\centering
\includegraphics[width=0.7\columnwidth]{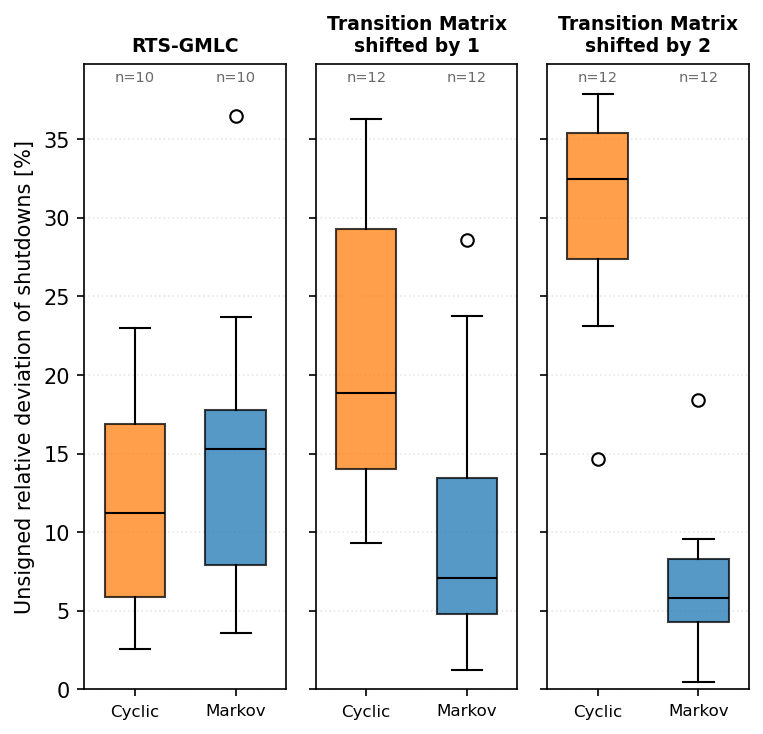}
\caption{Boxplots for the unsigned relative deviation of the number of shutdowns compared to the \Reference.}
\label{fig:vShutdownOperational}
\end{figure}

For the operational assessment, we fix the investments of all generators to the result of the \Reference solved as an investment problem, which provides the optimal investment\footnote{Note that the \Reference is then also resolved as an operational problem with fixed investments to obtain the correct work units for the comparison later. The \Reference is always solved for the full year, i.e., without RPs.}. Then, we solve the operational problem using the RP framework with different edge-handling methods for 3, 5, and 7 RPs, as well as demand-variability levels of 50\%, 70\%, 90\%, and 100\% (12 combinations in total\footnote{Note that the 90\% and 100\% demand-variability cases were excluded for RTS-GMLC with 3 RPs since the \Reference ran out-of-memory on our server. For the variants with shifted Transition Matrices, the \Reference for all 12 combinations could be solved.}). We define the \textit{operational deviation} of each method as the relative difference in the total number of shutdowns from the reference case, i.e., $(n_{\mathrm{method}}-n_{\mathrm{ref}})/n_{\mathrm{ref}}$, where $n$ denotes the respective count. Figure~\ref{fig:vShutdownOperational} shows the results for the shutdowns\footnote{Note that since the edges of RPs of \Markov and \Cyclic are connected, the number of shutdowns is equal to the number of start-ups.}. Results are shown for the original RTS-GMLC as well as two versions with shifted matrices. \NoEnf has been excluded from the plot, as its median deviation ranges from 68\% to 189\% and would distort the plot. Table~\ref{tab:rtsResults} shows the detailed results, reporting the median operational deviation, the median investment regret (defined in Subsection~\ref{subsubsec:InvestmentAssessment}), and the mean computational effort in work units relative to the full-year solve.

\begin{table}[!t]
\caption{Summary of the RTS-GMLC results.}
\label{tab:rtsResults}
\centering
\resizebox{\columnwidth}{!}{%
\begin{tabular}{@{}llrrrr@{}}
\hline
\multirow{2}{*}{\makecell{Transition\\Matrix}} & \multirow{2}{*}{\makecell{Edge-Handling\\Method}} & \multicolumn{2}{c}{Operational Assessment} & \multicolumn{2}{c}{Investment Assessment} \\ \cline{3-4}\cline{5-6}
 & & \makecell{Oper.\ dev.\\(median, \%)} & \makecell{Work units\\(mean, \% of Ref.)} & \makecell{Invest.\ regret\\(median, M\,EUR)} & \makecell{Work units\\(mean, \% of Ref.)} \\ \hline
\multirow{3}{*}{RTS-GMLC}    & No Enforcement      & $188.5$ & $0.5$ & $1.12$  & $0.5$ \\
                             & Cyclic Connection   & $-11.2$ & $0.3$ & $0.72$  & $0.1$ \\
                             & Markov Transition   & $15.3$  & $1.2$ & $-0.11$ & $0.1$ \\ \hline
\multirow{3}{*}{Shifted by 1} & No Enforcement      & $88.7$  & $1.6$ & $5.04$  & $1.4$ \\
                             & Cyclic Connection   & $-18.9$ & $1.3$ & $1.52$  & $0.2$ \\
                             & Markov Transition   & $7.1$   & $1.5$ & $-0.40$ & $0.2$ \\ \hline
\multirow{3}{*}{Shifted by 2} & No Enforcement      & $67.6$  & $1.4$ & $5.63$  & $1.0$ \\
                             & Cyclic Connection   & $-32.4$ & $1.1$ & $4.80$  & $0.3$ \\
                             & Markov Transition   & $5.8$   & $1.9$ & $0.00$  & $1.0$ \\ \hline
\end{tabular}}
\end{table}

Note that the \Cyclic consistently underestimates the number of shutdowns, whereas \NoEnf and the \Markov overestimate them\footnote{The plot shows the magnitude (absolute value) of the relative deviation to make comparison easier.}. For the original Transition Matrix, the \Markov performs comparably to the \Cyclic (a median deviation of $15\%$ versus $11\%$), as expected for a near-diagonal matrix. However, for the shifted Transition Matrices, the \Markov shows a significant improvement in the number of start-ups and shutdowns, reducing the magnitude of the median deviation from the reference case from $19\%$ and $32\%$ (\Cyclic) down to $7\%$ and $6\%$ (\Markov), respectively.

\begin{figure}[]
\centering
\includegraphics[width=1\columnwidth]{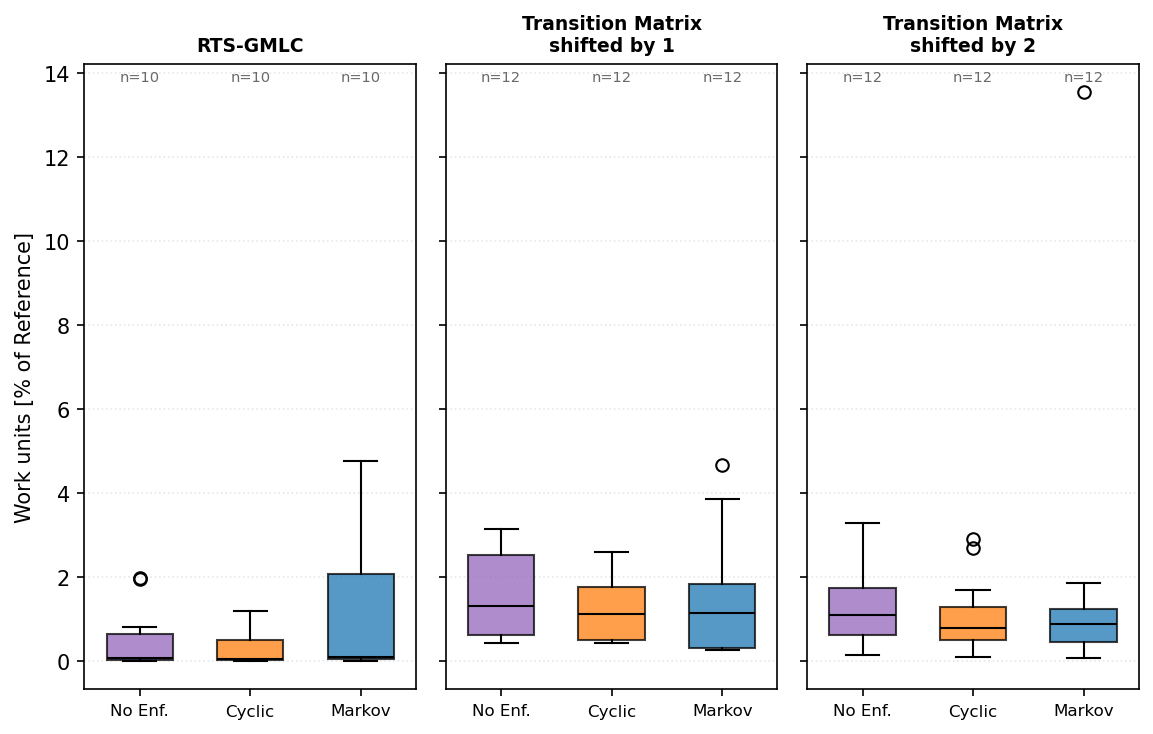}
\caption{Boxplots for the work units consumed relative to the \Reference in the operational assessment.}
\label{fig:WorkUnitsOperational}
\end{figure}

Regarding the work units consumed, Fig.~\ref{fig:WorkUnitsOperational} shows that the \Markov is either slower than the other two methods, or at best comparable. This is expected due to the additional interconnections between RPs (discussed in Subsection~\ref{subsec:Discussion}). The mean computational effort still stays below 2\% of that of the full model for all variants (also see Table~\ref{tab:rtsResults}), with only a single outlier reaching about 14\% (for the matrix shifted by 2). Relative to the \Cyclic, the \Markov adds at most $0.9$~percentage points of mean computational effort.

\subsubsection{Investment Assessment}\label{subsubsec:InvestmentAssessment}
To also assess the impact of the \Markov in investment decisions, we first solve the problem as an investment problem using RPs and the three different edge-handling methods. Then, we fix the resulting investment decisions in copies of the full model (with full temporal resolution) and solve them again as operational problems to obtain the objective function value. We define the \textit{investment regret} as the difference in total objective function value between such a full model with fixed (RP-based) investments and the \Reference, i.e., the full model solved as an investment problem without fixed investments. A regret of zero therefore means that the edge-handling method recovers investments that perform as well as those of a full chronological model. We again solve it for 3, 5, and 7 RPs, as well as demand-variability levels of 50\%, 70\%, 90\%, and 100\%\footnote{Again the two combinations already mentioned in Subsection~\ref{subsubsec:OperationalAssessment} had to be excluded due to an out-of-memory \Reference.}.

\begin{figure}[]
\centering
\includegraphics[width=1\columnwidth]{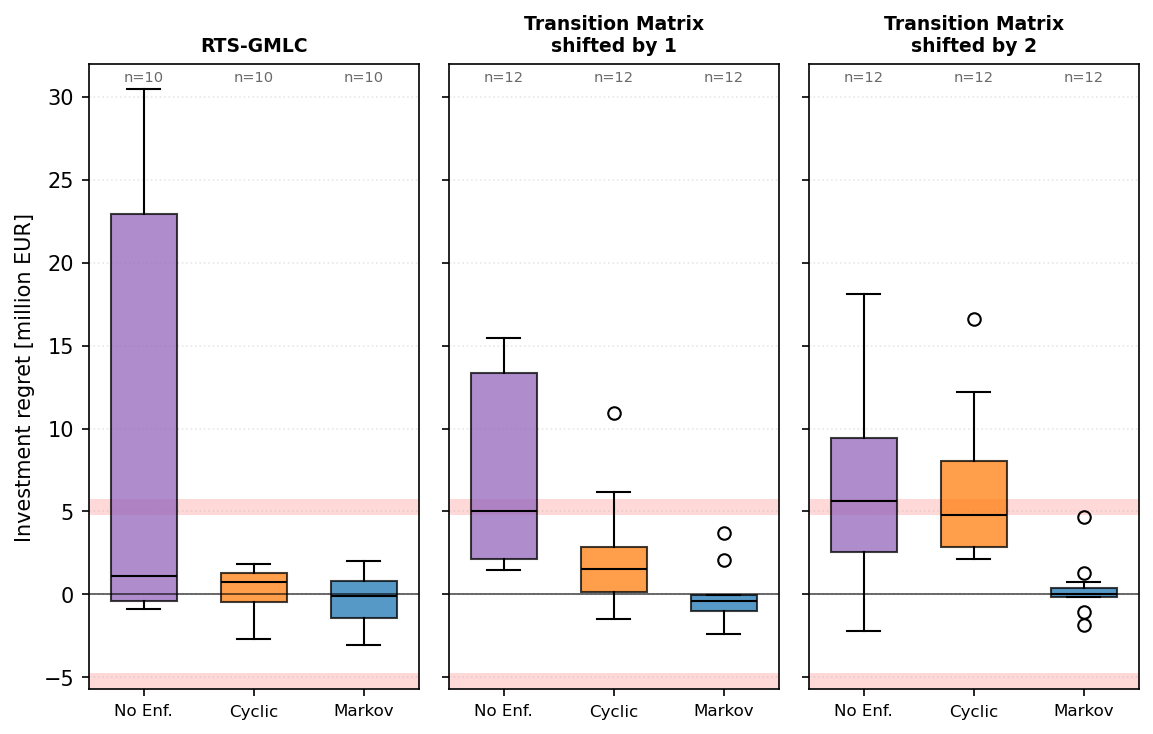}
\caption{Boxplots for the investment regret in million EUR. The red bands mark the MIP gap of $\pm0.3\%$ of the total cost.}
\label{fig:InvestmentAssessment}
\end{figure}

Figure~\ref{fig:InvestmentAssessment} shows this regret in million EUR\footnote{Slightly negative regret values are possible due to the MIP gap of 0.3\%, as the solver may stop at a solution that is not optimal.}. \NoEnf is the worst of the three, being up to about 30 million EUR more expensive. For the original Transition Matrix, \Markov and \Cyclic perform similarly, but for the shifted matrices, \Markov shows a significant improvement. Note, however, that the total objective is roughly 1900 million EUR across cases, meaning that any deviation below roughly 5 million EUR is below the MIP gap of 0.3\%. The MIP gap is shown as the red horizontal bands. The \Markov is the only method that consistently stays below this threshold, while the \Cyclic exceeds it for a substantial share of instances with the shifted matrices. See also Table~\ref{tab:rtsResults} for the exact median investment regret and average work units of the different edge-handling methods.

\begin{figure}[]
\centering
\includegraphics[width=1\columnwidth]{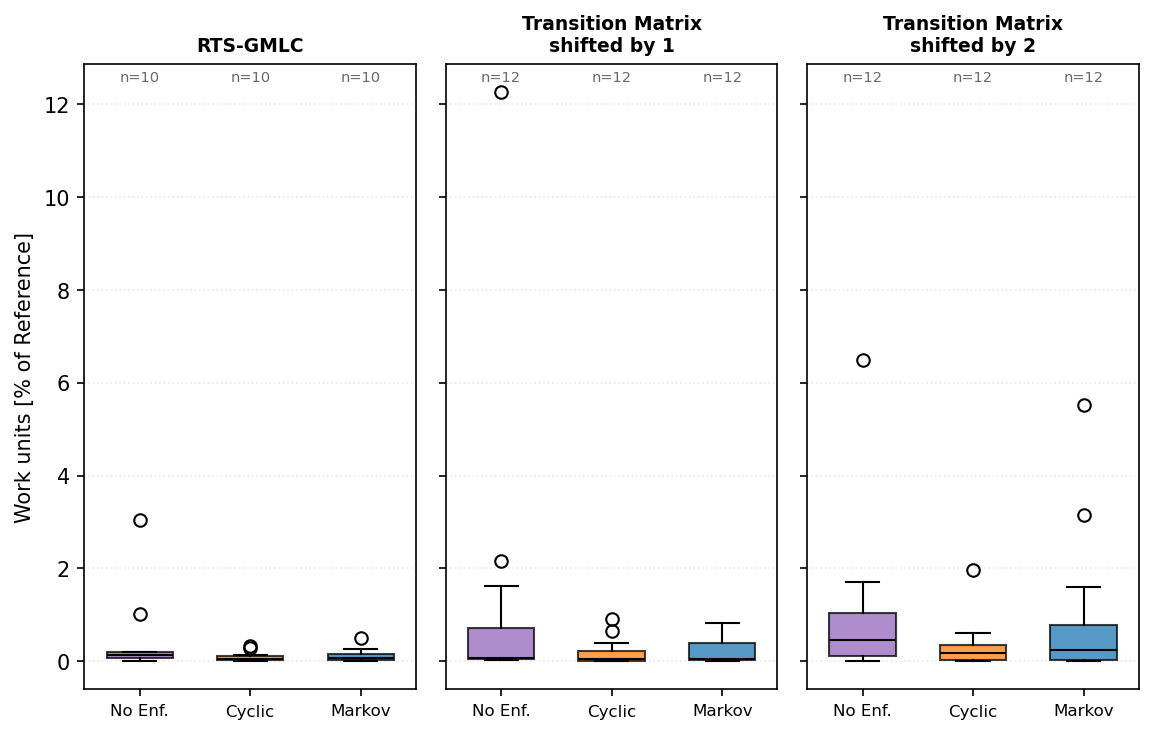}
\caption{Boxplots for the work units consumed relative to the \Reference in the investment assessment.}
\label{fig:InvestmentWorkUnits}
\end{figure}

Regarding the work units in the investment case, Fig.~\ref{fig:InvestmentWorkUnits} shows that the \Markov is on par with or slightly worse than the other two methods, and even its outliers never exceed 6\% of the work units consumed by the full model. Its mean effort differs from the \Cyclic by at most $0.7$~percentage points.

\subsection{Discussion}\label{subsec:Discussion}
The results show that the structure of the Transition Matrix has a significant impact on whether the \Markov can improve the results. Unsurprisingly, if the Transition Matrix is close to a diagonal matrix (which means that RPs tend to repeat themselves throughout the year), the \Cyclic is also very capable of approximating the results. For our case study on RTS-GMLC, the Transition Matrix is very close to a diagonal matrix for the tested number of RPs\footnote{Note that a higher number of RPs increases the probability of having a non-diagonal Transition Matrix.}. This leaves little room for improvement, with the \Markov and \Cyclic performing similarly. However, when we shift the Transition Matrix (to simulate a case study where RPs tend to iterate through each other), the \Markov shows a significant improvement over the \Cyclic\footnote{Note that preliminary results on case studies with random, dense Transition Matrices also show the same improved results when using the \Markov compared to the \Cyclic.}. We therefore recommend inspecting the Transition Matrix of a given dataset before committing to an edge-handling method. When its deviation from a diagonal matrix is large, the \Markov improves solution quality. Otherwise, the simpler \Cyclic formulation remains an appropriate choice. Constructing the matrix requires only the full chronological series that the a priori clustering already uses, so this inspection adds no data requirements.

The additional accuracy of the \Markov comes at a modest and bounded computational cost. Because every period boundary now depends on the boundaries of all other RPs through the Transition Matrix, the problem is no longer block-separable across periods. However, this cost is small, as the mean runtime stays below 2\% of the full-model runtime across all variants and exceeds that of the \Cyclic by at most 0.9 percentage points (Table~\ref{tab:rtsResults}).

One limitation concerns binary variables. Computing \emph{expected} boundary values can render commitment decisions fractional, which we try to mitigate by relaxing only the affected variables within the relevant up-/down-time window (Subsection~\ref{subsec:BinaryVariables}). In practice this affects fewer than $0.3\%$ of the binary variables: across all runs of the realistic case study, at most $0.28\%$ of the commitment, start-up, and shutdown flags take fractional values. Still, it means the returned schedule may be marginally non-integer and require a rounding or re-dispatch step before it is used directly as an operational plan. For the investment and aggregate-operational metrics studied here, this effect is negligible.

Finally, since the \Markov consistently overestimates start-ups and shutdowns, it can be considered the more conservative approach.

\section{Conclusions}\label{sec:Conclusion}
This paper introduces the \Markov formulation to enhance the chronological consistency of representative periods (RPs) by adjusting \textit{intra}-RP formulations. Comparative analyses against conventional formulations show that the \Markov formulation more accurately captures time-dependent system dynamics. Because periods that are critical for sizing the system may occur at RP boundaries, this improved edge representation propagates into the investment decisions as well, not just operations. For systems where RPs tend to appear in arbitrary succession and which have time-linked constraints active at RP boundaries, our approach has a significantly lower error, reducing the median operational deviation by up to $80\%$ (from about $32\%$ to $6\%$) in our realistic case study on the RTS-GMLC dataset. Notably, in the illustrative case, the \Markov is the only method that recovers the exact number and timing of all start-ups and shutdowns. These accuracy gains come at practically no extra cost, as the mean computational effort stays below $2\%$ of the full-model runtime, at most $0.9$ percentage points above the \Cyclic. In contrast, when such constraints are mostly active within individual RPs or if RPs tend to repeat themselves often, a simple \Cyclic formulation (where constraints connect the ends of each RP with the beginning of the same RP) remains an appropriate alternative.

Future work will focus on further improving scalability to very large systems. While the overhead is already small, the Transition Matrix couples all period boundaries; reducing it to its most significant entries (making it sparser) would cut the number of linking constraints and keep this overhead low as the number of RPs and the system size grow. The framework will also be tested on large, real-world datasets (e.g., on a model of the European network). Finally, the combination with spatial aggregation techniques is a natural direction of additional research, investigating trade-offs of temporal and spatial aggregation.

\section*{Acknowledgments}
The authors thank D. Cardona-Vasquez, S. Malacek, S. Karelly, and Y. Werner for their insightful discussions and constructive ideas.

\section*{Declaration of generative AI and AI-assisted technologies in the manuscript preparation process}
During the preparation of this work, the authors used Claude (Anthropic) to assist with writing the code for plots and tables, as well as for grammar checking and editorial refinement of the manuscript. The authors reviewed and edited the content as needed and take full responsibility for the content of the published article.

\printcredits

\section*{Declaration of competing interest}
The authors declare that they have no known competing financial interests or personal relationships that could have appeared to influence the work reported in this paper.

\section*{Data availability}
All data and code used for this study are available at \url{https://github.com/IEE-TUGraz/LEGO-Pyomo/tree/research/MarkovTransition}.






\bibliographystyle{unsrtnat}

\bibliography{references}



\end{document}